\documentclass[12pt]{amsart}

\begin{document}

\markboth{S. Yakubovich}{Index transforms with Kelvin functions}

\title{Index transforms with the squares of Kelvin functions} 

\author{S. Yakubovich}
\address{Department of Mathematics, Fac. Sciences of University of Porto,Rua do Campo Alegre,  687; 4169-007 Porto (Portugal)}




\maketitle

\begin{abstract}   New index transforms, involving squares of Kelvin functions,    are  investigated.   Mapping properties and inversion formulas are established for these transforms  in Lebesgue  spaces.  The results are applied to solve a  boundary value problem on the wedge for a  fourth   order partial differential equation.  
\end{abstract}

\vspace{6pt}


{\bf MS Classification}:  44A15, 33C10, 44A05\bigskip

\section{Introduction and preliminary results}

Let $ f(x),  \ g(\tau),\  x \in \mathbb{R}_+,\  \tau \in  \mathbb{R}_+$ be  complex-valued functions.  In this paper we will investigate mapping properties of the index transforms   \cite{yak} with the kernel  $ker^2_{2i\tau} \left( 2 (4x)^{1/4} \right) +
kei^2_{2i\tau} \left( 2 (4x)^{1/4} \right)$. Namely, we will consider the following integral operators
$$(F f ) (\tau) =  \int_0^\infty   \left[  ker^2_{2i\tau} \left( 2 (4x)^{1/4} \right) +
kei^2_{2i\tau} \left( 2 (4x)^{1/4} \right) \right] f(x)dx,  \eqno(1.1)$$
$$(G g) (x) =   \int_{0}^\infty   \left[  ker^2_{2i\tau} \left( 2 (4x)^{1/4} \right) +
kei^2_{2i\tau} \left( 2 (4x)^{1/4} \right) \right]   g(\tau) d\tau. \eqno(1.2)$$
Here  $i$ is the imaginary unit and $ker_\nu(x), kei_\nu(x)$ are Kelvin functions \cite{erd}, Vol. II, which are defined through the equality

$$ ker_\nu(x) \pm i kei_\nu(x) = e^{\mp  \nu \pi i/2} K_\nu \left(x e^{\pm \pi i/4}\right),\eqno(1.3)$$
where $K_\nu(z)$ is the modified Bessel function of the second kind \cite{erd}, Vol. II, which can be defined in terms of the integral

$$K_\nu(z)= \int_0^\infty e^{-z\cosh(y)}  \cosh(\nu y) dy,\  {\rm Re} z > 0.\eqno(1.4)$$ 
Hence the kernel of integral operators (1.1), (1.2) can be written in the form

$$  \mathcal{K}(x,\tau)  = ker^2_{2i\tau} \left( 2 (4x)^{1/4} \right) + kei^2_{2i\tau} \left( 2 (4x)^{1/4} \right) =  \left| K_{2i\tau} \left(  2 \left(4x e^{\pi i}\right)^{1/4} \right)\right|^2.$$
Another pair of Kelvin functions is defined by the equations 

$$  ber_\nu(x) \pm i bei_\nu(x) = e^{\pm  \nu \pi i}  J_\nu \left(x e^{\mp \pi i/4}\right) =  e^{\pm  \nu \pi i/2}  I_\nu \left(x e^{\pm \pi i/4}\right),\eqno(1.5)$$ 
where $J_\nu(z), I_\nu(z)$ are Bessel and modified Bessel functions of the first kind, respectively (cf. \cite{erd}, Vol. II ). 

Meanwhile, the Mellin-Barnes integral representation for the kernel $\mathcal{K}(x,\tau)$ is given according to relation (8.4.28.17) by the formula

$$ \mathcal{K}(x,\tau)  =  { 1 \over 16 \pi^{3/2} i } \int_{\gamma-i\infty}^{\gamma +i\infty}\Gamma(s+ i\tau)\Gamma(s-i\tau) $$

$$ \times \Gamma(s)  \Gamma(1/2+s) x^{-s} ds, \ x >0,\ \gamma >0,\eqno(1.6)$$
where $\Gamma(z)$ is Euler's gamma-function \cite{erd}, Vol. I  and $\gamma > 0$.    Integral (1.5) converges absolutely for any $\gamma > 0$ and uniformly with respect to $x \ge x_0 >0$ by virtue of the Stirling asymptotic formula for the gamma-function \cite{erd}, Vol. I, which yields for a fixed $\tau$

$$  \Gamma(s+ i\tau)\Gamma(s-i\tau)  \Gamma(s)  \Gamma(1/2+s) = O\left( e^{- 2\pi|s|} |s|^{ 4\gamma - 3/2}\right), \  |s| \to \infty.\eqno(1.7)$$
The Mellin-Barnes representation (1.6) can be used   to represent the kernel  for all $x >0,\ \tau \in \mathbb{R}$ as the Fourier cosine transform \cite{tit}.   Precisely, we have

{\bf Lemma 1}. {\it Let $x >0,  \tau \in \mathbb{R}$. Then the kernel  $\mathcal{K}(x,\tau)$  has the following representation in terms of the Fourier cosine transform }
$$   \mathcal{K}(x,\tau)  = \int_0^\infty  K_0\left( 4 x^{1/4}   \cosh^{1/2} (u)\right)   \cos(2u \tau) du.\eqno(1.8) $$

\begin{proof}     In fact,  appealing  to the reciprocal formulae via the Fourier cosine transform (cf. formula (1.104) in \cite{yak}) 
$$\int_0^\infty  \Gamma\left(s + i\tau\right)  \Gamma\left(s - i\tau \right)  \cos( \tau y) d\tau
= {\pi\over 2^{2s}}  {\Gamma(2s) \over \cosh^{2s}(y/2)},\ {\rm Re}\ s > 0,\eqno(1.9)$$
$$  \Gamma\left(s + i\tau \right)  \Gamma\left(s - i\tau\right)  
=   { \Gamma(2s)  \over 2^{2s-1}}  \int_0^\infty   {\cos(\tau y)  \over \cosh^{2s} (y/2)} \ dy,\eqno(1.10)$$ 
we replace  the gamma-product $ \Gamma\left(s + i\tau \right)  \Gamma\left(s - i\tau\right)$ in  (1.6)  by the right-hand side of  (1.10). Then, changing the order of integration via the absolute convergence, which can be justified using the Stirling asymptotic formula for the gamma-function, we employ the duplication formula  for the gamma-function \cite{erd}, Vol. I to derive 

$$   ker^2_{2i\tau} \left( 2 (4x)^{1/4} \right) + kei^2_{2i\tau} \left( 2 (4x)^{1/4} \right) =  { 1 \over 4 \pi i }  \int_0^\infty   \cos(\tau u)  \int_{\gamma-i\infty}^{\gamma +i\infty}  \left[ \Gamma(2s) \right]^2 $$

$$\times   \left(4 \sqrt x  \cosh (u/2)\right)^{- 2s  }  ds du.\eqno(1.11)$$
In the meantime, the inner integral with respect to $s$ in (1.11) is calculated in \cite{prud}, Vol. III, relation (8.4.23.1).  Hence we arrive at (1.7), completing the proof of Lemma 1. 

\end{proof}

Further,  recalling  the Mellin-Barnes representation (1.6) of the kernel  (1.8), we will derive an ordinary differential equation whose particular solution is  $ \mathcal{K} (x,\tau)$. Precisely,   it is given by 

{\bf Lemma 2}.    {\it   For each $\tau \in \mathbb{R}_+$ the function   $ \mathcal{K} (x,\tau)$ satisfies the following fourth  order differential equation with variable  coefficients}
$$x^3 {d^4  \mathcal{K} \over dx^4} +   {11\over 2} \  x^2 {d^3  \mathcal{K} \over dx^3} + x\left( {11\over 2} +\tau^2 \right) {d^2    \mathcal{K}  \over dx^2}  +  {1\over 2} \left( 1+\tau^2 \right) {d  \mathcal{K}  \over dx}  -  \mathcal{K}  = 0,\ x >0.\eqno(1.12)$$

\begin{proof}   Indeed,  the asymptotic behavior at infinity (1.7) of the integrand in (1.6) and the absolute and uniform convergence of the integral and its derivatives with respect to $x$ allow us to differentiate under the integral sign.     Hence, applying twice the differential operator $x{d\over dx} $ to both sides of (1.6) and employing the reduction formula for gamma-function \cite{erd}, Vol. I, we obtain  

$$ \left(x {d\over dx}\right)^2  \mathcal{K} (x,\tau) =   { 1 \over 16 \pi^{3/2} i } \int_{\gamma-i\infty}^{\gamma +i\infty} s^2 \Gamma(s+ i\tau)\Gamma(s-i\tau)  \Gamma(s)  \Gamma(1/2+s) x^{-s} ds$$

$$=  { 1 \over 16 \pi^{3/2} i } \int_{\gamma-i\infty}^{\gamma +i\infty}  \Gamma(s+1+ i\tau)\Gamma(s+1 -i\tau)  \Gamma(s)  \Gamma(1/2+s) x^{-s} ds - \tau^2  \mathcal{K} (x,\tau) $$

$$=  { 1 \over 16 \pi^{3/2} i } \int_{\gamma+1 -i\infty}^{\gamma +1 +i\infty}  \Gamma(s+ i\tau)\Gamma(s -i\tau)  \Gamma(s-1)  \Gamma(s- 1/2) x^{1 -s} ds - \tau^2  \mathcal{K} (x,\tau).$$
Making differentiation by $x$, it becomes 

$$ {d\over dx} \left(x {d\over dx}\right)^2  \mathcal{K} (x,\tau)  =  - { 1 \over 16 \pi^{3/2} i } \int_{\gamma+1 -i\infty}^{\gamma +1 +i\infty}  \Gamma(s+ i\tau)\Gamma(s -i\tau)  \Gamma(s)  \Gamma(s- 1/2) x^{-s} ds$$

$$ - \tau^2  {d\over dx} \mathcal{K} (x,\tau).$$
Now, multiplying by $x^{1/2}$ both sides of the obtained equality and differentiating again, we derive 

 $$\left( {d\over dx} x^{1/2}  {d\over dx} \right) \left(x {d\over dx}\right)^2  \mathcal{K} (x,\tau)  =   { 1 \over 16 \pi^{3/2} i } \int_{\gamma+1 -i\infty}^{\gamma +1 +i\infty}  \Gamma(s+ i\tau)\Gamma(s -i\tau)  \Gamma(s)  \Gamma(s +1/2)$$ 

 $$\times   x^{-s- 1/2} ds - \tau^2  \left( {d\over dx} x^{1/2}  {d\over dx} \right) \mathcal{K} (x,\tau).\eqno(1.13)$$
Meanwhile,  the contour in the integral on  the right-hand side of (1.12) can be moved to the left via Slater's theorem \cite{prud}, Vol. III. Hence, integrating again over the vertical line ${\rm Re}\  s = \gamma$ in the complex plane, we appeal to (1.5) to obtain the equality 

$$\left( {d\over dx} x^{1/2}  {d\over dx} \right) \left(\left(x {d\over dx}\right)^2 +\tau^2 \right)  \mathcal{K} (x,\tau) = 
{1\over \sqrt x}  \mathcal{K} (x,\tau).\eqno(1.14)$$
Finally, fulfilling the differentiation in (1.14), we end up with equation (1.12), completing the proof of Lemma 2. 

\end{proof} 

 In the sequel we will employ the Mellin  transform technique developed in \cite{yal} in order to  investigate the mapping properties of the index transforms  (1.1), (1.2).   Precisely, the Mellin transform is defined, for instance, in  $L_{\nu, p}(\mathbb{R}_+),\ 1 \le  p \le 2$ (see details in \cite{tit}) by the integral  
$$f^*(s)= \int_0^\infty f(x) x^{s-1} dx,\eqno(1.15)$$
 being convergent  in mean with respect to the norm in $L_q(\nu- i\infty, \nu + i\infty),\ \nu \in \mathbb{R}, \   q=p/(p-1)$.   Moreover, the  Parseval equality holds for $f \in L_{\nu, p}(\mathbb{R}_+),\  g \in L_{1-\nu, q}(\mathbb{R}_+)$
$$\int_0^\infty f(x) g(x) dx= {1\over 2\pi i} \int_{\nu- i\infty}^{\nu+i\infty} f^*(s) g^*(1-s) ds.\eqno(1.16)$$
The inverse Mellin transform is given accordingly
 $$f(x)= {1\over 2\pi i}  \int_{\nu- i\infty}^{\nu+i\infty} f^*(s)  x^{-s} ds,\eqno(1.17)$$
where the integral converges in mean with respect to the norm  in   $L_{\nu, p}(\mathbb{R}_+)$
$$||f||_{\nu,p} = \left( \int_0^\infty  |f(x)|^p x^{\nu p-1} dx\right)^{1/p}.\eqno(1.18)$$
In particular, letting $\nu= 1/p$ we get the usual space $L_p(\mathbb{R}_+; \  dx)$.   

The  Lebedev expansion theorem \cite{square}, which is related to the representation of the antiderivative in terms of the repeated integral with the product of the modified Bessel functions of the first and second kind \cite{erd}, Vol. II will play an important role to prove our main results.  Precisely,  we have 

{\bf Theorem 1}.  {\it Let $f \in L_{1/2,1} (0,1) \  \cap \  L_{3/2,1} (1,\infty) $. Then for all $ x >0$ 

$$\int_x^\infty f(y) dy= {4\over \pi^2} \int_0^\infty \tau \sinh(\pi\tau ) K_{i\tau}^2(x) d\tau$$

$$ \times \int_0^\infty f(y) \left[ I_{i\tau}(y) +  I_{- i\tau}(y)\right] K_{i\tau}(y) dy,\eqno(1.19)$$
 where   the inner integral by $y$ converges absolutely and the exterior integral by $\tau$ is understood as the Riemann improper integral. }

\section {Boundedness  and inversion properties of the index transform (1.1)}

Let $C_b (\mathbb{R}_+)$ be the space of  bounded continuous functions on $\mathbb{R}$.   

{\bf Theorem 2.}   {\it Let $ p \ge 1,\ q = p/(p-1),\ \nu <  1$. The index transform  $(1.1)$  is well-defined as a  bounded operator $F : L_{\nu,p} \left(\mathbb{R}_+\right) \to C_b (\mathbb{R}_+)$ and the following norm inequality takes place
$$||F f||_{C_b(\mathbb{R}_+)}  \le 4^{\nu- 3 + 1/q}  q^{ 4(\nu-1)} \  \Gamma^{1/q} ( 4(1-\nu) q) $$

$$\times  \ B(1-\nu, 1-\nu)  B(2 (1-\nu), 2(1-\nu))\  ||f||_{\nu,p},\eqno(2.1)$$
where 
$B(a,b)$ is the Euler beta-function \cite{erd}, Vol. I.   Moreover, it vanishes when $\tau \to +\infty$ .}

\begin{proof}  In fact,  appealing to Lemma 1, we prove  the following composition representation of the transform (1.1) in terms of the Fourier cosine and Meijer's  transforms \cite{yal}, namely,

$$(Ff)(\tau)= \int_0^\infty \cos(2u\tau)   \int_0^\infty  K_0\left( 4 x^{1/4}   \cosh^{1/2} (u)\right)  f(x) dx du.\eqno(2.2)$$
The interchange of the order of integration in (2.2) is allowed due to Fubini's theorem and the absolute convergence of the iterated integral. Indeed,  employing H{\"o}lder's, the generalized Minkowski inequalities and integral representation (1.4) of the modified Bessel function, we derive the estimates

$$ \int_0^\infty \left|\cos(2u\tau) \right|  \int_0^\infty  K_0\left( 4 x^{1/4}   \cosh^{1/2} (u)\right)  |f(x)|  dx du$$
 
$$\le  ||f||_{\nu,p}   \int_0^\infty \left( \int_0^\infty  K^q_0\left( 4 x^{1/4}   \cosh^{1/2} (u)\right)  x^{(1-\nu) q- 1} dx \right)^{1/q} du $$

$$= 4^{4(\nu-1) + 1/q}   ||f||_{\nu,p}   \int_0^\infty  {du \over \cosh^{2(1-\nu)} (u) } \left( \int_0^\infty  K^q_0\left( y\right)  y^{4(1-\nu) q- 1} dy \right)^{1/q}  $$

$$= 4^{3\nu- 4 + 1/q} \ B(1-\nu, 1-\nu)   ||f||_{\nu,p}    \left( \int_0^\infty \left(\int_0^\infty e^{- y \cosh (t)} dt\right)^q   y^{4(1-\nu) q- 1} dy \right)^{1/q}  $$

$$\le  4^{3\nu- 4 + 1/q} \ B(1-\nu, 1-\nu)   ||f||_{\nu,p}    \int_0^\infty  \left(\int_0^\infty e^{- y q\cosh (t)} y^{4(1-\nu) q- 1} dy \right)^{1/q   } dt  $$

$$=   4^{3\nu- 4 + 1/q}  q^{ 4(\nu-1)} \  \Gamma^{1/q} ( 4(1-\nu) q)  \ B(1-\nu, 1-\nu)   ||f||_{\nu,p}    \int_0^\infty  {dt \over \cosh^{4(1-\nu)}(t)}  $$

$$=   4^{\nu- 3 + 1/q}  q^{ 4(\nu-1)} \  \Gamma^{1/q} ( 4(1-\nu) q)  \ B(1-\nu, 1-\nu)  B(2 (1-\nu), 2(1-\nu))\  ||f||_{\nu,p} .$$
Hence,

$$ \sup_{\tau \in \mathbb{R}_+} \left| (F f)(\tau) \right| \le   4^{\nu- 3 + 1/q}  q^{ 4(\nu-1)} \  \Gamma^{1/q} ( 4(1-\nu) q) $$

$$\times  \ B(1-\nu, 1-\nu)  B(2 (1-\nu), 2(1-\nu))\  ||f||_{\nu,p}, $$
and we arrive at the estimate (2.1).  Moreover, as it follows from the latter estimate 

$$G(u)=  \int_0^\infty  K_0\left( 4 x^{1/4}   \cosh^{1/2} (u)\right)  f(x) dx \in L_1(\mathbb{R}_+).$$
Therefore appealing to the Riemann-Lebesgue lemma, we see from (2.2) that $(Ff)(\tau)$ vanishes at infinity.  Theorem 2 is proved. 

\end{proof} 

The inversion formula for the transform (1.1) is given by  

{\bf Theorem 3}. {\it Let $0 < \nu < 1/2$,  $f$ be locally integrable on $\mathbb{R}_+$, i.e. $f(x)\in L_{loc} (\mathbb{R}_+)$ and behave as  $f(x)= O\left(x^{-a}\right),\ x \to 0,\ f(x)= O\left(x^{-b}\right),\ x \to \infty,\  a < \nu < b,\ 3/2 > b > 5/4.$  Then if  its  Mellin transform \ $f^*(1) = f^{*\prime}(1) =0$ and $(Ff) (\tau) \in L_1\left(\mathbb{R}_+; \tau e^{\pi\tau} d\tau \right) $,  then the following inversion formula for the index transform $(1.1)$ holds for almost all $x >0$

$$    f(x)  =  -  {16 \over \pi } {d\over dx}   \int_{0}^\infty {\rm Im} \left[   \Gamma(1+2i\tau)  \left[ ber^2_{2i\tau} (2(4x)^{1/4} ) + bei^2_{2i\tau} (2 (4 x)^{1/4} ) \right]  \right] $$

$$\times  (Ff) (\tau) \tau  d\tau, \eqno(2.3)$$
where integral $(2.3)$ converges absolutely and ${\rm Im}$ denotes the imaginary part of a complex-valued function}. 

\begin{proof}  In fact,   recalling the   Mellin-Parseval equality (1.16), we  write (1.1) under conditions of the theorem in the form

$$(Ff) (\tau) =  { 1 \over 16 \pi^{3/2} i } \int_{\nu-i\infty}^{\nu  +i\infty}\Gamma(s+ i\tau)\Gamma(s-i\tau)  \Gamma(s)  \Gamma(1/2+s) f^*(1-s)  ds.\eqno(2.4)$$
Meanwhile, appealing to relation (8.4.23.23) in \cite{prud}, Vol. III, we have  

$${\sqrt \pi \over \cosh(\pi\tau)} K_{i\tau} \left(\sqrt x\right) \left[ I_{i\tau} \left(\sqrt x\right) + I_{-i\tau} \left(\sqrt x\right)\right]=
{1\over 2\pi i} \int_{\nu-i\infty}^{\nu  +i\infty}\Gamma(s+ i\tau)\Gamma(s-i\tau) $$

$$\times   \frac {  \Gamma(1/2-s)}{\Gamma(1-s)} x^{-s} ds, \ x> 0,\ 0< \nu < {1\over 2}.$$  
Hence via (2.4) and the Mellin-Parseval equality (1.6)

$$(Ff) (\tau) =  { 1 \over 8 \cosh(\pi\tau)} \int_0^\infty   K_{i\tau} \left(\sqrt x\right) \left[ I_{i\tau} \left(\sqrt x\right) + I_{-i\tau} \left(\sqrt x\right)\right] \varphi(x) dx,\eqno(2.5)$$
where

$$\varphi (x)=  {1\over 2\pi i} \int_{\nu-i\infty}^{\nu  +i\infty} \Gamma(1-s)\Gamma(3/2-s)\  {\Gamma(s) f^*(s) \over \Gamma(s-1/2) }  x^{-s} ds,\ x >0.\eqno(2.6) $$
A simple change of variables in (2.5) presumes

$$(Ff) (\tau) =  { 1 \over 4 \cosh(\pi\tau)} \int_0^\infty   K_{i\tau} \left( x\right) \left[ I_{i\tau} \left( x\right) + I_{-i\tau} \left( x\right)\right] \varphi(x^2) x  dx,\eqno(2.7)$$
In the meantime, from (2.6) and conditions of the theorem we get

$$\int_0^1 \sqrt x\ \left| \varphi (x^2) \right| dx = {1\over 2}  \int_0^1  \left| \varphi (x) \right| dx $$

$$\le  {1\over 4\pi (1-\nu) } \int_{\nu-i\infty}^{\nu +i\infty} \left| \Gamma(1-s)\Gamma(3/2-s)\  {\Gamma(s) f^*(s) \over \Gamma(s-1/2) }  ds\right| < \infty.$$
But the integrand in (2.6) is analytic in the vertical strip $a < \nu < b,\ 3/2 > b > 5/4$ because $f \in L_{\nu,1}(\mathbb{R}_+),\   a < \nu < b$ and $f^*(1)= f^{* \prime} (1)= 0$. Hence, shifting the contour in (2.6) to the right, we have the estimate
$$\int_1^\infty  x^{3/2} \ \left| \varphi (x^2) \right| dx = {1\over 2}  \int_1^\infty  x^{1/4}  \left| \varphi (x) \right| dx $$

$$\le  {1\over 4\pi (\nu- 5/4) } \int_{\nu -i\infty}^{\nu +i\infty} \left| \Gamma(1-s)\Gamma(3/2-s)\  {\Gamma(s) f^*(s) \over \Gamma(s-1/2) }  ds\right| < \infty,\ 5/4 < \nu < b.$$
Consequently, the function $x \varphi (x^2) $ satisfies conditions of the Lebedev expansion theorem (see Theorem 1), and we find

$$\int_x^\infty y \varphi (y^2) dy = {8 \over \pi^2} \int_0^\infty \tau \sinh(2 \pi\tau ) K_{i\tau}^2(x) (Ff)(\tau) d\tau .\eqno(2.8)$$
After  simple substitutions it becomes

$$\int_x^\infty  \varphi (y) dy = {16 \over \pi^2} \int_0^\infty \tau \sinh(2 \pi\tau ) K_{i\tau}^2(\sqrt x) (Ff) (\tau) d\tau,\ x >0.\eqno(2.9)$$
But from (2.6) and Fubini's theorem

$$\int_x^\infty  \varphi (y) dy = - {1\over 2\pi i} \int_{\nu-i\infty}^{\nu  +i\infty} \Gamma(1-s)\Gamma(3/2-s)\  {\Gamma(s) f^*(s) \over (1-s) \Gamma(s-1/2) }  x^{1-s} ds$$

$$ =   {1\over 2\pi i} \int_{1-\nu-i\infty}^{1-\nu  +i\infty} \Gamma(s)\Gamma(1/2+s)\  {\Gamma(-s) f^*(1-s) \over  \Gamma(1/2-s) }  x^{s} ds$$
Substituting the latter expression on the left-hand side of (2.9) and changing $x$ on $1/x$, we obtain 

$${1\over 2\pi i} \int_{1-\nu-i\infty}^{1-\nu  +i\infty} \Gamma(s)\Gamma(1/2+s)\  {\Gamma(-s) f^*(1-s) \over   \Gamma(1/2-s) }  x^{- s} ds$$

$$=  {16 \over \pi^2} \int_0^\infty \tau \sinh(2 \pi\tau ) K_{i\tau}^2\left({1\over \sqrt {x}} \right) (Ff) (\tau) d\tau.\eqno(2.10)$$
Meanwhile,  according to relation (8.4.23.28) in \cite{prud}, Vol. III the function $K_{i\tau}^2( 1/ \sqrt {x})$ has the following Mellin-Barnes integral representation

$$K_{i\tau}^2\left({1\over \sqrt x}\right) =  {1\over 4\pi^{1/2} i}  \int_{\gamma - i\infty}^{\gamma +i\infty}  \frac{ \Gamma( i\tau-s) \Gamma(-s- i\tau) \Gamma(-s)}{\Gamma(1/2 - s)}  x^{-s} ds,\  \gamma  < 0,\ x >0.\eqno(2.11)$$
Moreover, the Lebedev inequality (cf. \cite{yal}, p. 99) yields the estimate 

$$ K_{i\tau}^2\left({1\over \sqrt {x}} \right) \le {x^{1/4}  \over \sinh(\pi \tau)},\ x,\ \tau >0.\eqno(2.12)$$
Therefore, under assumed condition $(Ff)(\tau) \in L_1(\mathbb{R}_+;  \tau e^{\pi\tau} d\tau) $ we multiply both sides of (2.10) by $(t-x)^{-1/2} x^{-1/2},\ t > x $ and integrate with respect to $x$ over $(0,t)$.  The interchange of the order of integration is allowed in both sides of this equality by Fubini theorem with the use of (2.12) and Stirling asymptotic formula for the gamma-function.  Hence, employing (2.11) and  calculating  simple beta-integrals, we derive

$${1\over 2\pi i} \int_{1-\nu-i\infty}^{1-\nu  +i\infty} \Gamma(s)\Gamma(1/2+s)\  f^*(1-s) t^{- s} {ds \over s} $$

$$=  {4 \over \pi^{5/2} i} \int_0^\infty \tau \sinh(2 \pi\tau ) \int_{\gamma - i\infty}^{\gamma +i\infty}  \Gamma( i\tau-s) \Gamma(-s- i\tau) t^{-s} {ds\over s}  (Ff)(\tau) d\tau.\eqno(2.13)$$
Appealing to relation (8.4.23.2) in \cite{prud}, Vol. III, we have from (2.13) the equality

$${1\over 2\pi i} \int_{1-\nu-i\infty}^{1-\nu  +i\infty} \Gamma(s)\Gamma(1/2+s)\  f^*(1-s) t^{- s} {ds \over s} $$

$$= -  {16 \over \pi^{3/2} } \int_0^\infty \tau \sinh(2 \pi\tau ) \int_0^t K_{2i\tau} \left({2\over \sqrt x}\right) {dx\over x}  (Ff)(\tau) d\tau.\eqno(2.14)$$
The interchange of the order of integration on the right-hand side of (2.14)  is permitted due to the inequality (2.12) and the condition on $(Ff)(\tau) $. Indeed, we have the estimate

$$\int_0^\infty \tau \sinh(2 \pi\tau ) \int_0^t \left|K_{2i\tau} \left({2\over \sqrt x}\right) \right| {dx\over x}  |(Ff)(\tau) |d\tau$$

$$\le  \int_0^\infty {\tau \sinh(2 \pi\tau ) \over [ \sinh (2\pi\tau)]^{1/2} }|(Ff)(\tau) |d\tau   \int_0^t x^{-7/8} dx  < \infty.$$
Consequently, after differentiation with respect to $t$, we get from (2.14) 

$${1\over 2\pi i} \int_{1-\nu-i\infty}^{1-\nu  +i\infty} \Gamma(s)\Gamma(1/2+s)\  f^*(1-s) t^{- s}  ds $$

$$=   {16 \over \pi^{3/2} } \int_0^\infty \tau \sinh(2 \pi\tau )  K_{2i\tau} \left({2\over \sqrt t}\right) (Ff)(\tau) d\tau.\eqno(2.15)$$
Meanwhile, recalling the Mellin-Parseval equality (1.16) and the duplication formula for the gamma-function, we represent the modified Bessel function in the form

$$ K_{2i\tau} \left({2\over \sqrt x}\right) =    {1\over 4\pi i} \int_{\gamma - i\infty}^{\gamma +i\infty}   \Gamma(1+s)\Gamma(1/2+s) {\Gamma( i\tau-s) \Gamma(-s- i\tau)\over \Gamma(1+s)\Gamma(1/2+s)}  x^{-s} ds$$

$$ = {1\over 4 \sqrt \pi i} \int_{\gamma - i\infty}^{\gamma +i\infty}   \Gamma(1+2 s) {\Gamma( i\tau-s) \Gamma(-s- i\tau)\over \Gamma(1+s)\Gamma(1/2+s)}   (4 x)^{-s} ds$$

$$ =  {  \sqrt{\pi x}\over 4}   \int_0^\infty e^{- \sqrt {x y}} \  h(y,\tau) \sqrt y \ dy,\eqno(2.16)$$
where

$$h(x,\tau)  = {1\over 8\pi  i}  \int_{\mu - i\infty}^{\mu +i\infty}   {\Gamma( s-1 +i\tau) \Gamma(s-1- i\tau)\over \Gamma(2-s)\Gamma(3/2- s)} \   \left({x\over 4} \right)^{-s} ds,\ \mu  > 1. \eqno(2.17)$$
But the latter Mellin-Barnes integral converges absolutely for $1 < \mu  < 9/8$ and can be calculated by Slater's theorem  \cite{prud}, Vol. III in terms of the generalized hypergeometric functions ${}_0F_3 (a_1,a_2,a_3; z)$.  Thus we obtain  

$$ {1\over 2\pi i} \int_{\mu - i\infty}^{\mu +i\infty}   {\Gamma( s-1 +i\tau) \Gamma(s-1- i\tau)\over \Gamma(2-s)\Gamma(3/2- s)} \   \left({x\over 4} \right)^{-s} ds $$

$$=  \left({x\over 4} \right)^{i\tau -1} \sum_{n=0}^\infty {(-1)^n\over n!}   { \Gamma(-n - 2 i\tau)\over \Gamma(1+i\tau +n)\Gamma(1/2+i\tau +n)}   \left({x\over 4} \right)^n $$

$$+  \left({x\over 4} \right)^{- i\tau -1} \sum_{n=0}^\infty {(-1)^n\over n!}   { \Gamma(-n + 2 i\tau)\over \Gamma(1-i\tau +n)\Gamma(1/2-i\tau +n)}   \left({x\over 4} \right)^n $$

$$=  {2 \sqrt \pi \ x^{i\tau -1} \over  \tau \sinh(2\pi\tau) \Gamma(2i\tau) } \ {}_0F_3 \left( {1\over 2} + i\tau,\ 1+i\tau,\ 1+2i\tau ; \ {x\over 4} \right) $$  

$$+  {2 \sqrt \pi \ x^{- i\tau -1} \over \tau \sinh(2\pi\tau) \Gamma(- 2i\tau) } \ {}_0F_3 \left( {1\over 2} - i\tau,\ 1- i\tau,\ 1- 2i\tau ; \ {x\over 4} \right).\eqno(2.18) $$  
However,  relation (7.16.2.4) in \cite{erd}, Vol. III allows us to express the right-hand side of the latter equality in (2.18) in terms of Kelvin functions (1.5). Therefore

$${2 \sqrt \pi \ x^{i\tau -1} \over  \tau \sinh(2\pi\tau) \Gamma(2i\tau) } \ {}_0F_3 \left( {1\over 2} + i\tau,\ 1+i\tau,\ 1+2i\tau ; \ {x\over 4} \right) $$  

$$+  {2 \sqrt \pi \ x^{- i\tau -1} \over \tau \sinh(2\pi\tau) \Gamma(- 2i\tau) } \ {}_0F_3 \left( {1\over 2} - i\tau,\ 1- i\tau,\ 1- 2i\tau ; \ {x\over 4} \right)$$

$$= -  {8 \sqrt \pi \ \tau  \Gamma(2i\tau)  \over  x \sinh(2\pi\tau)  }  \left[ ber^2_{2i\tau} (2 x^{1/4} ) + bei^2_{2i\tau} (2 x^{1/4} ) \right] $$

$$ -  {8 \sqrt \pi \ \tau  \Gamma(- 2i\tau)  \over  x \sinh(2\pi\tau)  }  \left[ ber^2_{- 2i\tau} (2 x^{1/4} ) + bei^2_{- 2i\tau} (2 x^{1/4} ) \right]. \eqno(2.19)$$
Moreover, we derive the estimate 

$$|h(x,\tau) | \le  {\sqrt \pi \over 2 x } \left|  {x^{i\tau} \over  \tau \sinh(2\pi\tau) \Gamma(2i\tau) } \ {}_0F_3 \left( {1\over 2} + i\tau,\ 1+i\tau,\ 1+2i\tau ; \ {x\over 4} \right) \right. $$  

$$\left.  +   { x^{- i\tau} \over  \tau \sinh(2\pi\tau) \Gamma(-2i\tau) } \ {}_0F_3 \left( {1\over 2} - i\tau,\ 1-i\tau,\ 1-2i\tau ; \ {x\over 4} \right)  \right| $$  

$$\le  {\sqrt \pi \over 2 \tau \sinh(2\pi\tau) |\Gamma(2i\tau) | }  \left[ \sum_{n=0}^\infty {x^{n-1} \over 4^n [ n!] ^3 \  | (1/2 +i\tau)_n|  }  + \sum_{n=0}^\infty { x^{n-1} \over 4^n [ n!] ^3 \  | (1/2 - i\tau)_n|  } \right]$$

$$\le  { \sqrt \pi \over  \tau \sinh(2\pi\tau) |\Gamma(2i\tau) | }  \sum_{n=0}^\infty { x^{n-1} \over 2^n [ n!] ^3 \  (2n+1) !! } .$$
Returning to (2.16), we substitute its right-hand side into (2.15), having the estimate

$$\int_0^\infty \tau \sinh(2 \pi\tau ) \int_0^\infty e^{- \sqrt {x y}} \ | h(y,\tau) | \sqrt y \  |(Ff)(\tau)| dy  d\tau $$

$$\le   \sqrt \pi  \int_0^\infty  {   |(Ff)(\tau)|  d\tau \over   |\Gamma(2i\tau) | }  \int_0^\infty e^{- \sqrt {x y}} \   \sum_{n=0}^\infty { y^{n-1/2}  dy \over 2^n [ n!] ^3 \  (2n+1) !! } $$

$$=   2 \sqrt {{\pi\over x }}  \int_0^\infty  {   |(Ff)(\tau)|  \over   |\Gamma(2i\tau) | }  d\tau   \sum_{n=0}^\infty { x^{-n} \  (2n)! \over 2^n [ n!] ^3 \  (2n+1) !! }  < \infty,$$
which holds via Stirling asymptotic formulae for the gamma-function and the factorial, the condition  $(Ff)(\tau) \in L_1(\mathbb{R}_+;  \tau e^{\pi\tau} d\tau) $ and the continuity of $(Ff) (\tau)$ on $[0, +\infty)$ as it follows from (2.7) via the absolute and uniform convergence of the corresponding integral.  Hence by Fubini's theorem we justify the interchange of the order of integration. Therefore taking in mind  (2.18), (2.19) we get the equality

$${16 \over \pi^{3/2} } \int_0^\infty \tau \sinh(2 \pi\tau )  K_{2i\tau} \left({2\over \sqrt x}\right) (Ff) (\tau) d\tau$$

$$= -   16  {\sqrt{x \over \pi}}  \int_0^\infty e^{- \sqrt {x } y } \  dy  \int_0^\infty \tau^2 \left[   \Gamma(2i\tau)   \left[ ber^2_{2i\tau} (2 y^{1/2} ) + bei^2_{2i\tau} (2 y^{1/2} ) \right]  \right.$$

$$\left. +   \Gamma(- 2i\tau)    \left[ ber^2_{- 2i\tau} (2 y^{1/2} ) + bei^2_{- 2i\tau} (2 y^{1/2} )  \right] \right] (Ff) (\tau) d\tau.\eqno(2.20)$$
Treating similarly the left-hand side of (2.15), we find  

$${1\over 2\pi i} \int_{1-\nu-i\infty}^{1-\nu  +i\infty} \Gamma(s)\Gamma(1/2+s)\  f^*(1-s) x^{- s} ds $$

$$ = {1\over 2\sqrt \pi i} \int_{1-\nu-i\infty}^{1-\nu  +i\infty} \Gamma(1+2s)  f^*(1-s) \left(4x \right)^{- s} {ds\over s} $$

$$ = { \sqrt {\pi x}\over 8}  \int_0^\infty  e^{- \sqrt {xy} } \ {dy \over \sqrt y} \int_0^y f(u/4) du.$$
Consequently, combining with (2.20), we establish the following equality

$$  \int_0^\infty  e^{- \sqrt {x} y}   \int_0^{y^2/4}  f(u) du dy  = -  {16\over \pi}  \int_0^\infty e^{- \sqrt {x } y } \  dy  \int_{0}^\infty {\rm Im} \left[   \Gamma(1+2i\tau)  \right.$$

$$\times \left.  \left[ ber^2_{2i\tau} (2 y^{1/2} ) + bei^2_{2i\tau} (2 y^{1/2} ) \right]  \right]   (Ff) (\tau) \tau d\tau, x >0.\eqno(2.21)$$
Finally, the injectivity of the Laplace transform will be applied. Indeed,    equality (2.21) holds for all positive $t= \sqrt x$, and  therefore it is true for  a sequence of points $t_k= t_0+ k l,\  t_0, l >0,\ k= 1,2,\dots. $  Thus for all $x >0$

$$  \int_0^{x^2/4}  f(y) dy =  -  {16\over \pi}   \int_{0}^\infty   {\rm Im} \left[   \Gamma(1+2i\tau)   \left[ ber^2_{2i\tau} (2 x^{1/2} ) + bei^2_{2i\tau} (2 x^{1/2} ) \right]  \right]  (Ff) (\tau) \tau  d\tau.$$
Therefore,  making use simple substitutions  and differentiating the latter equality with respect to $x$, we arrive at the inversion formula (2.3) for almost all $ x>0$, completing the proof of Theorem 3. 

\end{proof}

\section{Index transform (1.2)} 

In this section we will study the  boundedness  properties and  prove an inversion formula for  the index transform  (1.2).   We begin with 

{\bf Theorem 4.}  {\it  Let $ \gamma  >0$ and  $ g(\tau)  \in L_1(\mathbb{R}_+;  d\tau)$. Then $x^{\gamma} (G g)(x)$ is bounded continuous on $\mathbb{R}_+$ and it holds  
$$\sup_{x >0}   \left| x^{\gamma}  (G g)(x) \right|   \le C_\gamma ||g ||_{L_1(\mathbb{R}_+; d\tau)},\eqno(3.1)$$
where

$$C_\gamma=  { 4^{- \gamma}  \over 8 \pi}  B(\gamma, \gamma ) \int_{\gamma-i\infty}^{\gamma +i\infty} \left| \Gamma^2 (s)  ds \right|.$$
Besides, if 
$$(Gg) (x) \in L_{\gamma,1}(\mathbb{R}_+),\  0 < \gamma < 1/2,$$
and its Mellin transform $|s|^{1/2-2\gamma} e^{\pi|s|/2} (Gg)^*(s) \in L_1( \gamma - i\infty,\  \gamma +i\infty)$, then for all $y >0$ }
$${1\over 2\pi i}   \int_{\gamma -i\infty}^{\gamma  +i\infty} \frac{\Gamma(1/2-s)}{\Gamma(s)\Gamma(1/2+s)} (G g)^*(s)  y^{ -s} ds$$

$$= { e^{y/2}   \over 8 }   \int_{0}^\infty   K_{i\tau} \left({y\over 2} \right) { g(\tau) \over \cosh (\pi\tau)}  d\tau. \eqno(3.2)$$

\begin{proof}   In fact,  substituting the right-hand side of (1.6) into (1.2), we estimate the obtained iterated integral, using (1.10) and the duplication formula for the gamma-function. As a result we get

$$ \left| x^{\gamma}  (G g)(x) \right|  \le   { 4^{-2\gamma}  \over 2  \pi } \int_{\gamma-i\infty}^{\gamma +i\infty} \left| \Gamma^2 (s)  ds \right| ||g||_{L_1(\mathbb{R}_+; d\tau)} \int_0^\infty {dy \over \cosh^{2\gamma} (y) } $$

$$= { 4^{- \gamma-1}  \over 2  \pi }  B(\gamma, \gamma ) \int_{\gamma-i\infty}^{\gamma +i\infty} \left| \Gamma^2 (s)  ds \right| ||g||_{L_1(\mathbb{R}_+; d\tau)}, $$
which yields (3.1).  Then taking the Mellin transform (1.15) of both sides  in (1.2) under the condition $(G g) (x) \in L_{\gamma,1}(\mathbb{R}_+),\ 0 < \gamma < 1/2$, we employ (1.6) and change the order of integration via Fubini's  theorem on the right-hand side of the obtained equality.  Hence we find 

$${  \Gamma(1/2-s) (Gg)^*(s) \over   \Gamma(s)  \Gamma(1/2+s) }=    {  \Gamma(1/2-s) \over 8 \sqrt \pi } \    \int_{0}^\infty    \Gamma(s+ i\tau)\Gamma(s-i\tau)   g(\tau) d\tau.\eqno(3.3)$$
Meanwhile, relation (8.4.23.5) in \cite{prud}, Vol. III says

$$e^{x/2} K_{i\tau} \left({x\over 2}\right) = {\cosh(\pi\tau)\over 2\pi^{3/2}i}   \int_{\gamma - i\infty}^{\gamma +i\infty}
 \Gamma(s+i\tau)\Gamma(s-i\tau) \Gamma(1/2-s) x^{-s} ds,\ 0< \gamma < {1\over 2}.$$
Hence an application of the inverse Mellin transform  (1.17) to both sides of (3.3) under conditions of the theorem drives us to (3.2).   The convergence of the integral with respect to $s$ in (3.2) is absolute under the condition $|s|^{1/2- 2\gamma} e^{\pi|s|/2} (G_\alpha g)^*(s) \in L_1( \gamma - i\infty,\  \gamma +i\infty)$, which  can be verified,  recalling the Stirling asymptotic formula for the gamma-function.  Theorem 4 is proved.

\end{proof} 

The inversion formula for the index transform (1.2) is given by

{\bf Theorem 5}.  {\it  Let  $g(z/i)$ be an even  analytic function in the strip $D= \left\{ z \in \mathbb{C}: \ |{\rm Re} z | \le \mu < 1/2\right\} ,\  g(0)= g^\prime(0) =0$ and $e^{-3 \pi | {\rm Im} z|/2} g(z/i)$ be  absolutely  integrable over any vertical line in  $D$.   If under conditions of Theorem 4 $(G g)^*(s)$ is analytic in $D$ and $|s|^{1/2- 2 \gamma} e^{\pi|s|/2} (G g)^*(s) \in L_1( \gamma - i\infty,\  \gamma +i\infty)$ over any vertical line ${\rm Re} s=\nu$ in the strip,  then for all  $x \in \mathbb{R}$ the  following inversion formula holds for the index transform (1.2)} 
$$ g(x)=    -  {16 x   \over  \pi  } \int_0^\infty  {\rm Im} \left[    \Gamma(1+2i x) \left[ ber^2_{2i x} (2 (4t)^{1/4} ) + bei^2_{2i x} (2 (4t)^{1/4} ) \right] \right] $$

$$\times  d  (Gg)(t).\eqno(3.4)$$

\begin{proof}    Indeed,  since the integrand on the left-hand side of (3.2) is analytic in the strip $|{\rm Re} s| \le \mu$ and absolutely integrable there over any vertical line, we shift the contour to the left, integrating over $(\nu-i\infty, \nu+i\infty), \nu < 0$.  Then multiplying both sides by $ e^{-y/2} K_{ix} \left({y/2} \right) y^{-1}$ and integrating  with respect to $y$ over $(0, \infty)$, we change the order of integration in the left-hand side of the obtained equality due to the absolute convergence of the  iterated integral.  Moreover,  appealing   to relation (8.4.23.3) in \cite{prud}, Vol. III, we calculate the inner integral  to find  the equality 
$$ {1\over 2\pi i}   \int_{\nu -i\infty}^{\nu  +i\infty} \frac {\Gamma(- s+ix)\Gamma(- s-ix) }{\Gamma(s) \Gamma(1/2 +s) } (G g)^*(s) ds$$
$$ =  {1\over 8 \sqrt \pi } \int_0^\infty  K_{ix} \left({y\over 2} \right) \int_{0}^\infty  K_{i\tau} \left({y\over 2} \right)  {g(\tau) \over \cosh (\pi\tau) } {d\tau dy\over y}. \eqno(3.5)$$
In the meantime the right-hand side of (3.5) can be written, employing  the  representation  of the Macdonald function in terms of the modified Bessel function of the first kind $I_z(y)$ \cite{erd}, Vol. II 

$${2 i \over \pi}  \sinh(\pi\tau)  K_{i\tau} \left({y\over 2} \right) =   I_{ -i\tau} \left({y\over 2} \right) -  I_{i\tau} \left({y\over 2} \right),$$
the substitution $z=i\tau$ and the property  $g(- z/i)=  g(z/i)$.  Hence it becomes

$${1\over 8 \sqrt \pi } \int_0^\infty  K_{ix} \left({y\over 2} \right)  \int_{0}^\infty  K_{i\tau} \left({y\over 2} \right)  {g(\tau) \over \cosh (\pi\tau) } {d\tau dy\over y}$$

$$ = -  {\sqrt \pi \over 8 i   } \int_0^\infty  K_{ix} \left({y\over 2} \right)  \int_{-i\infty}^{i\infty}   I_{ z} \left({y\over 2} \right)  g\left({z\over i}\right) {dz  dy\over y \sin(2\pi z) }.\eqno(3.6)$$
On the other hand, according to our assumption $g(z/i)$ is analytic in the vertical  strip $0\le  {\rm Re}  z \le \mu,\  0 < \mu < 1/2$,  and $e^{- 3\pi |{\rm Im} z|/2} g(z/i)$ is absolutely integrable  in the strip.  Hence,  appealing to the inequality for the modified Bessel   function of the first  kind  (see \cite{yal}, p. 93)
 $$|I_z(y)| \le I_{  {\rm Re} z} (y) \  e^{\pi |{\rm Im} z|/2},\   0< {\rm Re} z \le  \mu,$$
one can move the contour to the right in the latter integral in (3.6). Then 

$$-  {\sqrt \pi \over 8 i   }  \int_0^\infty  K_{ix} \left({y\over 2} \right)  \int_{-i\infty}^{i\infty}   I_{ z} \left({y\over 2} \right)  g\left({z\over i}\right) {dz  dy\over y \sin(2\pi z) }$$

$$= -  {\sqrt \pi \over 8 i   }  \int_0^\infty  K_{ix} \left({y\over 2} \right)  \int_{\mu -i\infty}^{\mu+ i\infty}   I_{ z} \left({y\over 2} \right)  g\left({z\over i}\right) {dz  dy\over y \sin(2\pi z) }.\eqno(3.7)$$
Hence one can  interchange the order of integration in the right-hand side of (3.7) due to the absolute and uniform convergence.  Then using  the value of the integral (see relation (2.16.28.3) in \cite{prud}, Vol. II)
$$\int_0^\infty K_{ix}(y) I_z(y) {dy\over y} = {1\over x^2 + z^2},$$ 
we find  

$$ -  {\sqrt \pi \over 8 i   }  \int_0^\infty  K_{ix} \left({y\over 2} \right)  \int_{\mu -i\infty}^{\mu+ i\infty}   I_{ z} \left({y\over 2} \right)  g\left({z\over i}\right) {dz  dy\over y \sin(2\pi z) }$$

$$=   -  {\sqrt \pi \over 8 i   }  \int_{\mu -i\infty}^{\mu+ i\infty}  {g\left(z/ i\right) \over (x^2+ z^2) \sin(2\pi z) } dz $$

$$= -  {\sqrt \pi \over 8 i   } \left( \int_{-\mu +i\infty}^{- \mu- i\infty}   +   \int_{\mu -i\infty}^{ \mu +  i\infty}   \right)  {  g(z/i) \  dz \over z (z-ix) \sin(2\pi z)}.\eqno(3.8)$$
Thus we are ready to apply the Cauchy formula in the right-hand side of the latter equality in (3.8) under conditions of the theorem.  Hence 

$$ -  {\sqrt \pi \over 8 i   }  \int_0^\infty  K_{ix} \left({y\over 2} \right) \int_{\mu -i\infty}^{\mu+ i\infty}   I_{ z} \left({y\over 2} \right)  g\left({z\over i}\right) {dz  dy\over y\sin(2\pi z)}$$

$$=  {\pi^{3/2}  \over 8   }   \ { g(x) \over  x \sinh(2\pi x) } ,\quad x \in \mathbb{R} \backslash \{0\}.\eqno(3.9)$$
Combining with (3.5), we establish the equality 

$$      {1\over 2\pi i} \int_{\nu -i\infty}^{\nu  +i\infty} \frac {\Gamma(- s+ ix)\Gamma(- s-ix) }{\Gamma(s) \Gamma(1/2+s) } (G g)^*(s) ds$$

$$=   {\pi^{3/2}  \over 8   }   \ { g(x) \over  x \sinh(2\pi x) } ,\quad x \in \mathbb{R} \backslash \{0\}.\eqno(3.10)$$
The left-hand side of (3.10) can be treated, in turn, appealing to the Parseval equality (1.16), conditions of the theorem and the value of the integral (2.18) (see (2.19)).   Therefore we derive

$$      {1\over 2\pi i} \int_{\nu -i\infty}^{\nu  +i\infty} \frac {\Gamma(- s+ ix)\Gamma(- s-ix) }{\Gamma(s) \Gamma(1/2+s) } (G g)^*(s) ds$$

$$  =    {1\over 2\pi i} \int_{1-\nu -i\infty}^{1-\nu  +i\infty} \frac {\Gamma(s-1+ ix)\Gamma( s-1-ix) }{\Gamma(2-s) \Gamma(3/2-s) }  (1-s) (G g)^*(1-s) ds$$

$$= -  {2 x \sqrt \pi \over  \sinh(2\pi x) }\int_0^\infty  \left[    \Gamma(2i x) \left[ ber^2_{2i x} (2 (4t)^{1/4} ) + bei^2_{2i x} (2 (4t)^{1/4} ) \right] \right.$$

$$\left. +    \Gamma(- 2i x)  \left[ ber^2_{- 2i x} (2 (4t)^{1/4} ) + bei^2_{- 2i x} (2 (4t)^{1/4} ) \right] \right]   d (Gg)(t).$$
Consequently, recalling (3.10), we establish the inversion formula  for the index transform (1.2)

$$ g(x)=    -  {8 x   \over  \pi i } \int_0^\infty  \left[    \Gamma(1+2i x) \left[ ber^2_{2i x} (2 (4t)^{1/4} ) + bei^2_{2i x} (2 (4t)^{1/4} ) \right] \right.$$

$$\left. -     \Gamma(1- 2i x)  \left[ ber^2_{- 2i x} (2 (4t)^{1/4} ) + bei^2_{- 2i x} (2 (4t)^{1/4} ) \right] \right]  d (Gg)(t),$$
which is equivalent to (3.4). Theorem 4 is proved.

 \end{proof}

\section{Boundary   value problem}

In this section we will employ  the index transform (1.2) to investigate  the  solvability  of a boundary  value  problem  for the following fourth    order partial differential  equation, involving the Laplacian

$$  \left(  x {\partial^2  \over \partial x^2}  + y  {\partial^2   \over  \partial x \partial y} \right)^2 u   +  \left(  y {\partial^2  \over \partial y^2}  + x  {\partial^2   \over  \partial x \partial y} \right)^2 u $$ 

$$+ \left( {9\over 2} \left(  x {\partial  \over \partial x}  + y  {\partial   \over  \partial y} \right)+ 3 \right)  \Delta u -  {u \over r }  = 0, \quad   (x,y) \in  \mathbb{R}^2 \backslash \{0\},\eqno(4.1)$$ 
where $r= \sqrt{x^2+y^2},\ \Delta = {\partial^2 \over \partial x^2} +  {\partial^2 \over \partial y ^2}$ is the Laplacian in $\mathbb{R}^2$.   In fact, writing  (4.1) in polar coordinates $(r,\theta)$, we end up with the equation  

$$r^3 {\partial^4  u \over \partial r^4} +  r {\partial^4  u \over \partial r^2 \partial \theta^2 \ }+   {11\over 2} \  r^2 {\partial^3  u \over \partial r^3}  + {1\over 2} \ {\partial^3  u \over \partial r \partial \theta^2} +  {11\over 2}  r  {\partial^2 u  \over \partial r^2}     +  {1\over 2}  {\partial u \over \partial r}  -  u = 0.\eqno(4.2)$$

{\bf Lemma 3.} {\it  Let $g(\tau)  \in L_1\left(\mathbb{R}_+;  ( \tau^2 +1) d\tau\right),\  \beta \in (0, 2\pi)$. Then  the function
$$u(r,\theta)=    \int_{0}^\infty    \left[  ker^2_{2i\tau} \left( 2 (4r)^{1/4} \right) +
kei^2_{2i\tau} \left( 2 (4r)^{1/4} \right) \right] \  {\sinh(\theta \tau)\over \sinh(\beta\tau) }   g(\tau) d\tau,\eqno(4.3)$$
 satisfies   the partial  differential  equation $(4.2)$ on the wedge  $(r,\theta): r   >0, \  0\le \theta \le  \beta$, vanishing at infinity.}

\begin{proof} The proof  is straightforward by  substitution (4.3) into (4.2) and the use of  (1.12).  The necessary  differentiation  with respect to $r$ and $\theta$ under the integral sign is allowed via the absolute and uniform convergence, which can be verified  analogously to estimates in Theorem 4  under the  condition $g \in L_1\left(\mathbb{R};  (\tau^2 +1) d\tau\right)$.   Finally,  the condition $ u(r,\theta) \to 0,\ r \to \infty$  is due to the integral representation (1.8) of the kernel  in (4.3). 
\end{proof}

Finally  we will formulate the boundary  value problem for equation (4.2) and give its solution.

{\bf Theorem 6.} {\it Let  $g(x)$ be given by formula $(3.4)$ and its transform $(G g) (t)\equiv G (t)$ satisfies conditions of Theorem 5.  Then  $u (r,\theta),\   r >0,  \  0\le \theta \le \beta$ by formula $(4.3)$  will be a solution  of the boundary  value problem for the partial differential  equation $(4.2)$ subject to the conditions}
$$u(r,0) = 0,\quad\quad   u(r,\beta) =  G (r).$$

\bigskip
\noindent{{\bf Funding}}
\bigskip

\noindent The work was partially supported by CMUP (UID/MAT/00144/2013),  which is funded by FCT (Portugal) with national (MEC),   European structural funds through the programs FEDER  under the partnership agreement PT2020, and Project STRIDE - NORTE-01-0145-FEDER- 000033, funded by ERDF - NORTE 2020. 

\bibliographystyle{amsplain}

\begin{thebibliography}{10}

\bibitem{yak}   Yakubovich S.  Index transforms.   Singapore:  World Scientific Publishing Company; 1996.

\bibitem{erd}    Erd\'elyi A,  Magnus W,   Oberhettinger  F,   Tricomi FG.  Higher transcendental functions. Vols. I,  II. New  York: McGraw-Hill;  1953.

\bibitem{prud}  Prudnikov AP,  Brychkov  YuA,  Marichev OI. Integrals and series:  Vol. I: Elementary functions. New York:  Gordon and Breach;   1986;   Vol. II:  Special functions. New York: Gordon and Breach;  1986;   Vol. III:  More special functions. New York:   Gordon and Breach; 1990.

\bibitem{yal}   Yakubovich S,  Luchko Yu.  The hypergeometric approach to integral transforms and convolutions, Mathematics and its applications.  Vol. 287.  Dordrecht:  Kluwer Academic Publishers Group; 1994.

\bibitem {tit}  Titchmarsh EC.   An introduction to the theory of Fourier integrals.   New York:  Chelsea; 1986.

\bibitem{square}  Lebedev  NN.   On an integral representation of an arbitrary function in terms of squares of Macdonald functions with imaginary index.   Sibirsk. Mat. Zh. 1962; 3: 213- 222 (in Russian).





\end{thebibliography}

\end{document}